\documentclass[]{article}

\usepackage{type1cm}        
\usepackage{makeidx}         
\usepackage{graphicx}        
\usepackage{multicol}        
\usepackage[bottom]{footmisc}
\usepackage{todonotes}
\usepackage{hyperref}
\usepackage[a4paper, total={15cm, 25cm}]{geometry}

\usepackage{newtxtext}       %
\usepackage[varvw]{newtxmath}       

\usepackage{subcaption}

\makeindex             

\usepackage{bm}

\usepackage{xcolor}

\newcommand{\pf}{\varphi}

\newcommand{\strain}[1][u]{{\bm \varepsilon}(\bm #1)}
\newcommand{\Div}{\nabla\cdot}

\usepackage{color,soul}

\usepackage{tikz}
\usepackage{pgfplots,pgfplotstable}
\pgfplotsset{compat=newest}

\definecolor{blau}{RGB}{0 144 188}
\definecolor{newblue1}{RGB}{0 144 188}
\definecolor{newblue2}{RGB}{197 216 227}
\definecolor{newgreen1}{RGB}{0 144 118}
\definecolor{newgreen2}{RGB}{197 222 215}
\definecolor{neworange1}{RGB}{255 137 0}
\definecolor{neworange2}{RGB}{255 205 105}
\definecolor{newred1}{RGB}{254 54 41}
\definecolor{newred2}{RGB}{141 26  18}
\definecolor{newpurple1}{RGB}{196 19 252}
\definecolor{newpurple2}{RGB}{93 14 117}

\usepackage{authblk}

\title{Sequential solution strategies for the Cahn-Hilliard-Biot model}
\author[1]{Erlend Storvik \footnote{Corresponding author: erlend.storvik@hvl.no}}
\author[2]{Cedric Riethm\"uller}
\author[3]{Jakub W. Both} 
\author[3]{Florin A. Radu}

\date{}

\affil[1]{Department of Computer science, Electrical engineering and Mathematical sciences, Western Norway University of Applied sciences, Norway} 
\affil[2]{Institute of Applied Analysis and Numerical Simulation, University of Stuttgart, Germany}
\affil[3]{Centre for Modeling of Coupled Subsurface Dynamics, Department of Mathematics, University of Bergen, Norway}

\begin{document}

\maketitle
\subsection*{Abstract}
This paper presents a study of solution strategies for the Cahn-Hilliard-Biot equations, a complex mathematical model for understanding flow in deformable porous media with changing solid phases. Solving the Cahn-Hilliard-Biot system poses significant challenges due to its coupled, nonlinear and non-convex nature.
We explore various solution algorithms, comparing monolithic and splitting strategies, focusing on both their computational efficiency and robustness.

\section{Introduction}
The Cahn-Hilliard-Biot model introduced in \cite{Storvik2022} is a three-way coupled system that features the interplay of solid phase separation, fluid dynamics, and elastic deformations in porous media. Recently, well-posedness of the model was discussed in \cite{Riethmüller2023}, whereas the analysis for a similar model was performed in \cite{Fritz2023}. The present paper explores numerical solution strategies where we delve into both monolithic and splitting methods to solve the discretized equations. We investigate both computational efficiency and robustness with respect to material parameters. The monolithic approach, applying Newton's method in an implicit time-discretization, showcases efficiency in convergence but faces challenges in robustness. In contrast, the splitting strategy sequentially addresses the Cahn-Hilliard subsystem, the elasticity subsystem, and the flow subsystem, offering flexibility and adaptability in managing the unique dynamics of each component.

The Cahn-Hilliard-Biot system presents particular challenges due to its coupled, nonlinear, and non-convex characteristics. In this paper, we focus on handling the couplings in the system by comparing a decoupling method with a monolithic solution approach. This system's nonlinearities are predominantly manifested in its coupling terms. Therefore, by implementing a three-way decoupling strategy, we essentially introduce a method that serves not only to decouple but also to partially linearize the system. 

Comprehensive research specifically dedicated to solution strategies for the Cahn-Hilliard-Biot equations appears to be relatively limited. This work intends to contribute to this area, acknowledging that there is much yet to be explored and understood. Relevant to this context, there have been some studies on related models. These include \cite{Storvik2023} on time-discretization for the Cahn-Larch\'e equations, a sequential accelerated approach for phase-field fracture modeling discussed in \cite{Storvik2021}, strategies for optimal stabilization in decoupling the Biot equations in \cite{Storvik2019}, and considerations of solution strategies for splitting schemes in soft material poromechanics as explored in \cite{both2022}.

The paper is structured as follows. In Section 2 we introduce the mathematical model, in Section 3 we discuss the discretization in time and space and present the splitting schemes. Section 4 is devoted to illustrative numerical results. The paper ends with a concluding section. 

\section{The Cahn-Hilliard-Biot model}
Let $\Omega \subset \mathbb{R}^d$ be a domain with a Lipschitz continuous boundary, $d \in \{2,3\}$ be the spatial dimension and $T_\mathrm{f} > 0$ be the final time. We consider the problem: Find $(\pf, \mu, \bm u, p, {\bm q})$ for $({\bf x},t)\in \Omega \times [0,T_\mathrm{f}]$, with $\varphi$ being the phase-field variable, $\mu$ the chemical potential, $\bm u$ the displacement, $p$ the pore pressure, and $\bm q$ the fluid flux, such that
\begin{subequations}
\begin{align}
\partial_t \varphi - \Div (m \nabla \mu) &= R, \label{eq:ch1}\\
\mu + \gamma\left(\ell\Delta\varphi-\frac{1}{\ell}\Psi'(\varphi)\right) - \delta_\varphi\mathcal{E}_\mathrm{e}(\varphi, \bm u) - \delta_\varphi \mathcal{E}_\mathrm{f}(\varphi, \bm u, p) &= 0, \label{eq:ch2}\\
 -\nabla\cdot\big(\mathbb{C}(\varphi)\left(\bm\varepsilon\left(\bm u\right)-\mathcal{T}(\varphi)\right)-\alpha(\varphi)p\bm I\big) &= {\bm f}, \label{eq:elasticity}\\ 
\partial_t\left(\frac{p}{M(\varphi)} + \alpha(\varphi)\nabla\cdot \bm u\right) + \Div {\bm q} &= S_\mathrm{f}, \label{eq:flow}\\
\kappa(\pf)^{-1} {\bm q} +\nabla p &= 0,\label{eq:darcyflow}
\end{align}
\label{eq:model}%
\end{subequations}
where 
$$\delta_\varphi\mathcal{E}_\mathrm{e}(\varphi, \bm u) = -T'\!(\varphi)\!:\!\mathbb{C}(\varphi)\big(\bm\varepsilon\left(\bm u\right)-\mathcal{T}(\varphi)\big)+\frac{1}{2} \big(\strain - \mathcal{T}(\varphi)\big)\!:\!\mathbb{C}'\!(\varphi)\big(\strain - \mathcal{T}(\varphi)\big),$$
and 
$$\delta_\varphi\mathcal{E}_{\mathrm{f}}(\varphi,\bm u, p) = \frac{M'(\varphi)p^2}{2M(\varphi)^2}-p\alpha'(\varphi)\nabla\cdot u$$
accompanied with appropriate initial and boundary conditions. 
Here, $\strain := \frac{1}{2}\left(\nabla \bm u + \nabla\bm u^\top\right)$ is the linearized strain tensor, $\mathcal{T}(\pf):= \xi(\pf - \bar{\pf})\bm I$ accounts for swelling effects, $m$ is the chemical mobility, $\mathbb{C}(\varphi)$ is the stiffness tensor, $\alpha(\pf)$ is the Biot-Willis coupling coefficient, $M(\pf)$ is the compressibility coefficient, $\gamma > 0$ denotes the surface tension, $\ell > 0$ is a small parameter associated with the width of the regularization layer, and $\Psi(\varphi)$ is a double-well potential penalizing deviations from the pure phases, having minimal values at $\pf = 1$ and $\pf = 0$. In this paper it is given as $\Psi(\varphi) := \varphi^2\left(1-\varphi\right)^2.$ We note that this set of equations consists of a Cahn-Hilliard like subsystem \eqref{eq:ch1}--\eqref{eq:ch2}, linearized elasticity with momentum balance \eqref{eq:elasticity} and single-phase Darcy flow with mass balance \eqref{eq:flow}--\eqref{eq:darcyflow}.

\section{Numerical solution strategies}
In this paper we will investigate two different solution strategies; the monolithic Newton-based method, described in Section~\ref{sec:monolithic}, and a decoupling algorithm, described in Section~\ref{sec:splitting}. 

\subsection{Time-discretization: Implicit with partly convex-concave split}\label{sec:time_disc}
Both of the solution strategies that are presented here, utilize a semi-implicit time-discretization, with all terms being evaluated implicitly, except for the double-well potential $\Psi(\varphi).$ The double-well is split into a combination of a convex term and an expansive term as follows, $\Psi(\varphi) = \Psi_c(\varphi)- \Psi_e(\varphi)$, where both $\Psi_c(\varphi)$ and $\Psi_e(\varphi)$ are convex. Then, the convex term $\Psi_c(\varphi)$ is evaluated implicitly, while the expansive term $\Psi_e(\varphi)$ is evaluated explicitly, taking inspiration from the approach in \cite{eyre1998}. This split is not unique, and in this paper we use
\begin{equation}
      \Psi(\varphi) = \varphi^2(\varphi-1)^2 = \left(\left(\varphi-\frac{1}{2}\right)^4+\frac{1}{16}\right)-\frac{1}{2}\left(\varphi-\frac{1}{2}\right)^2= \Psi_c(\varphi)-\Psi_e(\varphi).
\end{equation}

\subsection{Spatial discretization}\label{sec:spatial_disc}
We apply a mixed finite element discretization of the system \eqref{eq:ch1}--\eqref{eq:darcyflow}. In particular, the Cahn-Hilliard subsystem \eqref{eq:ch1}--\eqref{eq:ch2}, is discretized by first order Lagrange elements for both the phase-field $\varphi$ and the chemical potential $\mu$, the elasticity subsystem \eqref{eq:elasticity} by vectorial first order Lagrange elements, and the flow subsystem \eqref{eq:flow}--\eqref{eq:darcyflow} using piecewise constant elements for the pore pressure $p$ and lowest-order Raviart-Thomas elements for the flux $\bm q$, ensuring local mass conservation. In the following we use the notation $\mathcal{V}_h^\mathrm{ch}$, $\mathcal{V}_h^{\bm u}$, $\mathcal{V}_h^{p}$, and $\mathcal{V}_h^{\bm q}$ for the test and trial spaces for $\varphi$ and $\mu$, $\bm u$, $p$, and $\bm q$, respectively.

\subsection{Monolithic solution strategy}\label{sec:monolithic}

We examine a monolithic solution strategy, wherein Newton's method is employed to solve the aforementioned discretized system of equations, see Sections~\ref{sec:time_disc}--\ref{sec:spatial_disc}. This approach involves applying Newton's method to the entire system of equations collectively, rather than decomposing it into smaller segments. While this holistic method is theoretically appealing for its ability to maintain the interrelated dynamics of the equations, it's important to acknowledge that Newton's method can suffer from a lack of robustness, especially in scenarios with complex nonlinearities, such as this. Here, the nonlinear iterations are initialized with the solution to the previous time step.

\subsection{Splitting schemes}\label{sec:splitting}
In contrast to the monolithic approach, this paper also explores a splitting strategy for the Cahn-Hilliard-Biot equations, which offers an alternative route by decomposing the complex system into more manageable subproblems. This strategy involves sequentially solving the Cahn-Hilliard subsystem \eqref{eq:ch1}--\eqref{eq:ch2}, the elasticity subsystem \eqref{eq:elasticity}, and the flow subsystem \eqref{eq:flow}--\eqref{eq:darcyflow} in an iterative loop. By tackling these subsystems individually, the strategy can utilize specialized numerical methods tailored to the specific challenges of each subproblem. Therefore, the sequential approach allows for greater flexibility and adaptability in handling the unique properties of each equation. In the implementation used in this paper, the Cahn-Hilliard subsystem is solved using the Newton method in each sequential iteration, whereas both the elasticity and flow subproblems become linear when decoupled and only need to be solved once. A similar strategy is applied, properly presented, and identified with alternating minimization in \cite{Storvik2023}. We note, that it is often beneficial, and necessary, to add a stabilizing term when applying a splitting scheme, as carefully discussed in \cite{Storvik2019}. We have, however, not applied any kind of stabilization to the decoupling scheme here.

After discretization, using index $n$ to indicate the time step, and $i$ as the iteration index, the iterative decoupling scheme reads: Given $\varphi_h^{n-1}, \varphi_h^{n,i-1}\in \mathcal{V}_h^{\mathrm{ch}}$, $\bm u_h^{n,i-1}, \bm u_h^{n-1}\in \mathcal{V}_h^{\bm u}$, and $p_h^{n,i-1}, p_h^{n-1}\in \mathcal{V}_h^{p}$ find $\varphi_h^{n,i}, \mu_h^{n,i}\in \mathcal{V}_h^{\mathrm{ch}}$, $\bm u_h^{n,i}\in \mathcal{V}_h^{\bm u},$ $p_h^{n,i}\in \mathcal{V}_h^{p}$ and $\bm q_h^{n,i}\in \mathcal{V}_h^{\bm q}$ such that
\begin{subequations}
\begin{align}
\left(\varphi^{n,i}_h-\varphi_h^{n-1},\eta_h^\varphi\right) + m\left( \nabla \mu_h^{n,i}, \nabla \eta_h^\varphi\right)- \tau \left(R,\eta_h^\varphi\right) &= 0, \label{eq:ch1split}\\
\left(\mu^{n,i},\eta_h^\mu\right) - \gamma\left(\ell\left(\nabla\varphi_h^{n,i},\nabla\eta^\mu_h\right)+\frac{1}{\ell}\left(\Psi_c'\left(\varphi_h^{n,i}\right)-\Psi_e'\left(\varphi_h^{n-1}\right),\eta_h^\mu\right)\right)\nonumber\\ - \left(\delta_\varphi\mathcal{E}_\mathrm{e}\left(\varphi_h^{n,i}, \bm u_h^{n,i-1}\right) - \delta_\varphi \mathcal{E}_\mathrm{f}\left(\varphi_h^{n,i}, \bm u_h^{n,i-1}, p_h^{n,i-1}\right),\eta_h^\mu\right) &= 0, \label{eq:ch2split}\\
\left(\mathbb{C}\left(\varphi_h^{n,i}\right)\left(\bm\varepsilon\left(\bm u_h^{n,i}\right)-\mathcal{T}\left(\varphi_h^{n,i}\right)\right),\bm\varepsilon\left(\bm \eta_h^{\bm u}\right)\right)\nonumber\\-\left(\alpha\left(\varphi_h^{n,i}\right)p_h^{n,i-1}, \bm \nabla\cdot \bm \eta_h^{\bm u}\right)-\left({\bm f},\bm \eta_h^{\bm u}\right) &= 0, \label{eq:elasticitysplit}\\ 
\left(\frac{p^{n,i}_h}{M\left(\varphi_h^{n,i}\right)} + \alpha\left(\varphi_h^{n,i}\right)\nabla\cdot \bm u_h^{n,i} -\left(\frac{p^{n-1}_h}{M\left(\varphi_h^{n-1}\right)} + \alpha\left(\varphi_h^{n-1}\right)\nabla\cdot \bm u_h^{n-1}\right), \eta_h^p\right)\nonumber\\+ \tau\left(\Div {\bm q_h^{n,i}}, \eta_h^p\right)- \tau\left(S_\mathrm{f},\eta_h^p\right) &= 0, \label{eq:flowsplit}\\
\left(\kappa\left(\pf_h^{n,i}\right)^{-1} {\bm q_h^{n,i}},\bm\eta_h^{\bm q}\right) -\left( p_h^{n,i}, \nabla\cdot \bm \eta_h^{\bm q}\right) &= 0,\label{eq:darcyflowsplit}
\end{align}
\label{eq:model}%
\end{subequations}
for all $\eta_h^\varphi, \eta_h^\mu\in \mathcal{V}_h^\mathrm{ch},$ $\bm \eta_h^{\bm u}\in \mathcal{V}_h^{\bm u}$, $\eta_h^p\in \mathcal{V}_h^p$ and $\bm \eta_h^{\bm q}\in \mathcal{V}_h^{\bm q}$. The system is then solved by forward substitution where first the Cahn-Hilliard subsystem \eqref{eq:ch1split}--\eqref{eq:ch2split} is solved, using Newton's method to resolve the nonlinearities, to obtain $\varphi_h^{n,i}$ and $\mu_h^{n,i}$. Then, the, now linear, elasticity equation \eqref{eq:elasticitysplit} is solved to obtain $\bm u_h^{n,i}$, before the, now also linear, flow subsystem \eqref{eq:flowsplit}--\eqref{eq:darcyflowsplit} is solved to obtain $p_h^{n,i}$ and $\bm q_h^{n,i}$. The method proceeds iteratively until a stopping criterion is reached, and the iterations are initialized by using the solutions at the previous time step.

\section{Numerical experiments}
In this section, we present a series of numerical experiments designed to evaluate the efficacy of the proposed solution strategies for the Cahn-Hilliard-Biot equations. These experiments are crafted to be suitable for understanding how well both of the presented schemes, monolithic and splitting, are performing dependent on changing material parameters. The nonlinear material parameters that depend on the solid phase are written in terms of the interpolation function
\begin{equation}
    \label{eq:interpolation}
    \pi(\varphi) = \begin{cases}
        0, &\mathrm{for}\quad \varphi<0,\\
        \varphi^2\left(3-2\varphi\right), &\mathrm{for}\quad \varphi\in\left[0,1\right],\\
        1, &\mathrm{for}\quad \varphi>1,\\
    \end{cases}
\end{equation}
as $\zeta(\varphi) = \zeta_0 + \pi(\varphi)\left(\zeta_1-\zeta_0\right),$ with $\zeta$ being a placeholder for $\mathbb{C},$ $M$, $\alpha$, and $\kappa$. The specific material, discretization and linearization parameters are presented in Table~\ref{tab:1}, with the stiffness tensors given as 
\begin{equation}\label{eq:stiffness}
    \mathbb{C}_{0} = 
    \begin{pmatrix}
        100 & 20  & 0\\
        20  & 100 & 0\\
        0   & 0   & 100
    \end{pmatrix}, \quad \mathrm{and}\quad 
    \mathbb{C}_{1} = 
    \begin{pmatrix}
        1   & 0.1 & 0\\
        0.1 & 1   & 0\\
        0   & 0   & 1
    \end{pmatrix},
\end{equation}
in Voigt notation.

\begin{table}[h]
\centering
\begin{tabular}{c|c|c|c}
Parameter name                    & Symbol                         & Value         &  Unit \\
\hline
Surface tension                   & $\gamma$                       & $5$           & $\left[F\right]$\\
Regularization parameter          & $\ell$                         & $2.0$e$-2$    & --\\
Mobility                          & $m$                            & $1.0$         & $\left[\frac{L^4}{FT}\right]$  \\
Swelling parameter                & $\xi$                          & $0.5$         & --\\
Reference phase-field             & $\bar{\varphi}$              & $0.5$         & --\\
Stiffness tensors                 & $\mathbb{C}_0$, $\mathbb{C}_1$ & -             & $\left[\frac{F}{L^2}\right]$\\
Compressibilities                 & $M_0$, $M_1$                   & $1.0$, $0.1$  & $\left[\frac{F}{L^2}\right]$\\
Permeabilities                    & $\kappa_0$, $\kappa_1$         & $1.0$, $0.1$  & $\left[\frac{L^4}{FT}\right]$ \\
Biot-Willis coupling coefficients & $\alpha_0$, $\alpha_1$         & $1.0$, $0.5$  & -- \\
Time step size                    & $\tau$                         & $1$e$-5$      & $[T]$\\
Mesh size                         & $h$                            & $\sqrt{2}/65$ & $[L]$\\
Max nonlinear iteration           & max\textunderscore iter        & $100$         & -- \\
Tolerance                         & tol                            & $1$e$-6$      & -- \\
\end{tabular}
\caption{Simulation parameters. The stiffness tensors are given in \eqref{eq:stiffness} and the interpolation function in \eqref{eq:interpolation}. Note that the value for $\gamma$ is varying in Section~\ref{sec:surf}, and the value for $\xi$ is varying in Section~\ref{sec:swell}. 
Here, $L$, $T$ and $F$ denote the units of length, time and force respectively. 
}
\label{tab:1}
\end{table}

\begin{figure}
\centering
    \begin{subfigure}{0.19\textwidth}
        \includegraphics[width = \textwidth]{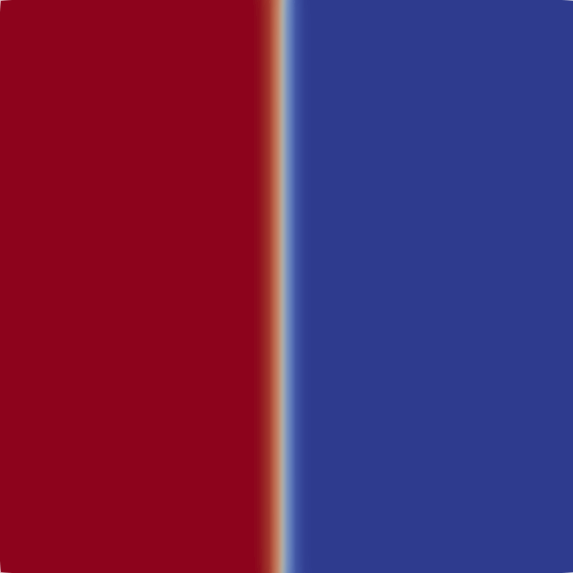}
        \caption{$t = 0$}
        \label{fig:initial}
    \end{subfigure}
    \hfill
    \begin{subfigure}{0.19\textwidth}
        \includegraphics[width = \textwidth]{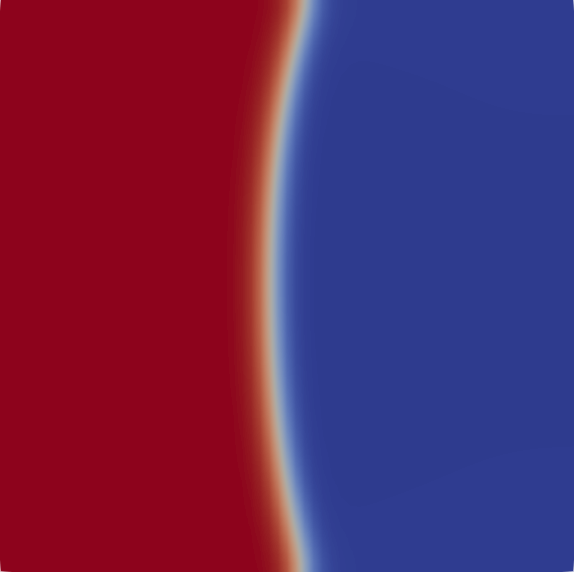}
        \caption{$\xi = 0.2$}
        \label{fig:end01}
    \end{subfigure}
    \hfill
    \begin{subfigure}{0.19\textwidth}
        \includegraphics[width = \textwidth]{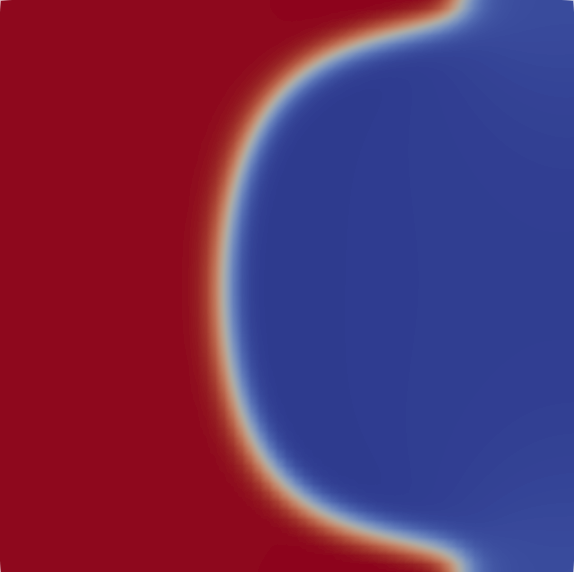}
        \caption{$\xi = 0.5$}
        \label{fig:end025}
    \end{subfigure}
    \hfill
    \begin{subfigure}{0.19\textwidth}
        \includegraphics[width = \textwidth]{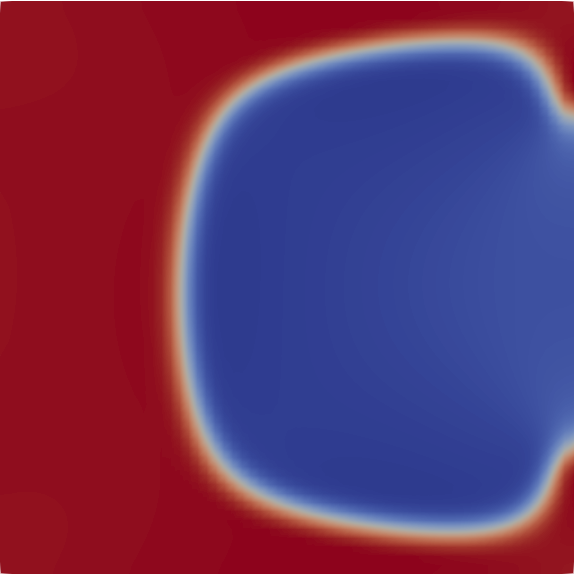}
        \caption{$\xi = 1.0$}
        \label{fig:end05}
    \end{subfigure}
    \hfill
    \begin{subfigure}{0.19\textwidth}
        \includegraphics[width = \textwidth]{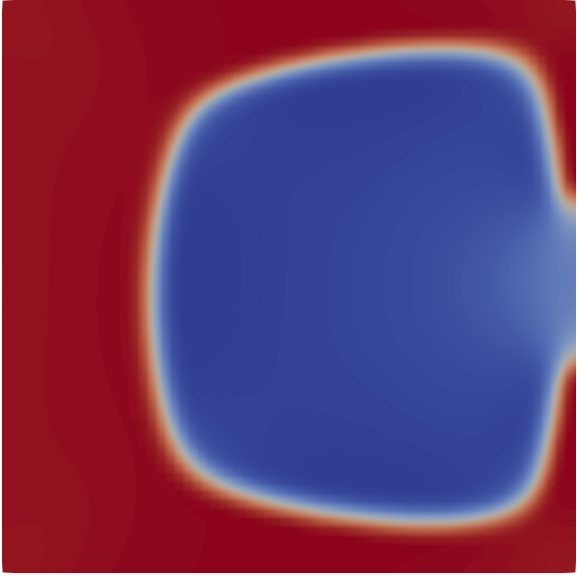}
        \caption{$\xi = 1.5$}
        \label{fig:end075}
    \end{subfigure}
    \hfill
\label{fig:CHB}
       \caption{Plot of phase-field solutions $\varphi$ from numerical study. The phase-field was initialized with the distribution from figure (a) for all cases. The other plots (b)--(e) show the phase-field at the final time $T_\mathrm{f} = 0.003$ for different values of the swelling parameter $\xi$.}
\end{figure}

For all of the experiments the spatial domain $\Omega$ is the unit square. We apply homogeneous Dirichlet boundary conditions to the elasticity and flow subproblems, and homogeneous Neumann conditions to the Cahn-Hilliard subproblem. The source terms $R$, $\bm f$ and $S_\mathrm{f}$ are all set to zero. The simulation is initialized with a phase-field that takes value zero in one half of the domain and one in the other half, see Figure~\ref{fig:initial}, which eventually drives the dynamics. The pressure and displacement variables are initialized as $p({\bf x},0)=0$ and $\bm u({\bf x},0) = \bf 0.$

\subsection{Implementation}
The implementation is done using the FEniCS computing platform \cite{alnaes2015fenics, logg2012automated}, and the source code for this project is available at \href{https://github.com/EStorvik/chb}{https://github.com/EStorvik/chb}.

\subsection{Dependence on surface tension parameter}\label{sec:surf}
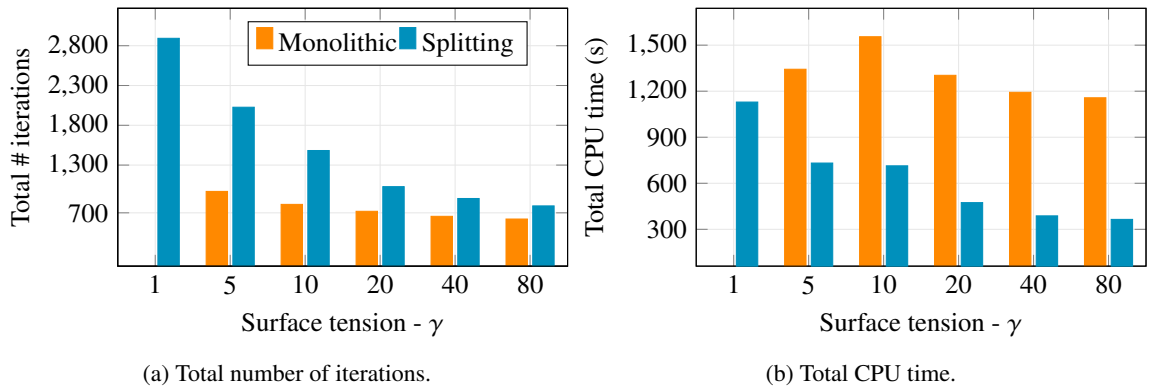
\begin{figure}
\begin{subfigure}{0.5\textwidth}
\begin{tikzpicture}
   \pgfplotsset{ybar stacked, ymin=300.000, ymax=3000.000, width=\textwidth, enlargelimits=0.1,
   grid=both, grid style={line width=.2pt, draw=gray!20}, 
   xmin=0, xmax = 5,
 }
 \begin{axis}
 [
 ylabel={Total \# iterations}, 
 xticklabels={1, 5, 10, 20, 40, 80},
 xlabel = {Surface tension - $\gamma$},
 ytick={700, 1300, 1800, 2300, 2800}, 
 height=5cm,
 xtick=data,
 bar width=8pt,
 bar shift=-5pt, 
 height=5cm,
 scaled y ticks = false
 ]  
 
 \pgfplotsinvokeforeach {1}{
     \addplot[color=neworange1, fill=neworange1] table [x expr=\coordindex, y index=#1] {\datatableGammaMonoIter};}
     
\end{axis} 

 \begin{axis}
 [
 bar width=8pt,
 bar shift=5pt, 
 xtick=\empty,
 ytick=\empty, 
 height=5cm
 ]  
 
 \pgfplotsinvokeforeach {1}{
     \addplot[color=newblue1, fill=newblue1] table [x expr=\coordindex, y index=#1] {\datatableGammaSplitIter};}
     
\end{axis} 

 \begin{axis}
 [
 bar width=0pt,
 bar shift=10pt, 
 legend columns=2,
 legend style={at={(0.6,0.95)}, anchor=north}, 
 xtick=\empty,
 legend cell align=left,
 legend style={
     /tikz/column 2/.style={
         column sep=-1.3pt,
         row sep=-20pt 
     },
     /tikz/column 6/.style={ 
         column sep=-1.3pt,
         row sep=-20pt 
     }},
 ytick=\empty,
 height=5cm
 ]  
 \pgfplotsinvokeforeach {1}{
     \addplot[color=neworange1, fill=neworange1] table [x expr=\coordindex, y index=#1] {\datatableEmpty};}
 \pgfplotsinvokeforeach {1}{
     \addplot[color=newblue1, fill=newblue1] table [x expr=\coordindex, y index=#1] {\datatableEmpty};}

 \legend{
 Monolithic, 
 Splitting, 
 }
 \end{axis}
 \end{tikzpicture}
 \caption{Total number of iterations.}
 \label{fig:gammaiterations}
\end{subfigure}
\begin{subfigure}{0.5\textwidth}
  \begin{tikzpicture}
   \pgfplotsset{ybar stacked, ymin=200.000, ymax=1600.000, width=\textwidth, enlargelimits=0.1,
   grid=both, grid style={line width=.2pt, draw=gray!20}, 
   xmin=0, xmax = 5,
 }
 
 \begin{axis}
 [
 ylabel={Total CPU time (s)}, 
 xticklabels={1, 5, 10, 20, 40, 80},
 xlabel = {Surface tension - $\gamma$},
 ytick={300, 600, 900, 1200, 1500}, 
 height=5cm,
 xtick=data,
 bar width=8pt,
 bar shift=-5pt, 
 height=5cm,
 scaled y ticks = false
 ]  
 
 \pgfplotsinvokeforeach {1}{
     \addplot[color=neworange1, fill=neworange1] table [x expr=\coordindex, y index=#1] {\datatableGammaMonoTime};}
     
\end{axis} 

 \begin{axis}
 [
 bar width=8pt,
 bar shift=5pt, 
 xtick=\empty,
 ytick=\empty, 
 height=5cm
 ]  
 
 \pgfplotsinvokeforeach {1}{
     \addplot[color=newblue1, fill=newblue1] table [x expr=\coordindex, y index=#1] {\datatableGammaSplitTime};}
     
\end{axis} 


 \end{tikzpicture}
  \caption{Total CPU time.}
 \label{fig:gammatime}
 \end{subfigure}
 \caption{The total number of iterations and CPU time to complete the simulation for different values of the surface tension parameter $\gamma$. The material parameters from Table~\ref{tab:1} are applied with swelling parameter $\xi = 0.5$.}
 \label{fig:gamma}
\end{figure}

The monolithic and splitting schemes proposed in this study were evaluated across a range of surface tension parameter $\gamma$ values. Figure~\ref{fig:gamma} displays the total iterations and overall time (in seconds) required for the complete simulation. It is observed that both schemes exhibit a similar trend: higher surface tension values simplify the resolution of nonlinear problems. This simplification is attributed to the increased prominence of the Cahn-Hilliard subproblem, which becomes more amenable to linearization when subjected to convex-concave time discretization, resulting in a well-conditioned system for linear methods.

Note that the monolithic solver failed to converge for the lowest value of the surface tension parameter $\gamma$, whereas the splitting method succeeded. Moreover, although the splitting method required more iterations for convergence, it was faster in terms of total CPU time spent, indicating its efficiency. It's important to note that the implementation did not employ specialized techniques for assembling or solving the linear systems. For straightforward implementations, such as the one undertaken in this study, the decoupled solver proved to be significantly quicker and more reliable. While the linear systems in the monolithic solver could potentially benefit from customized linear solvers and preconditioners, the comparison remains equitable as both the monolithic and splitting methods utilized the same direct methods and straightforward implementations.

\pgfplotstableread{
0 0
1 0
80 0
}\datatableEmpty

\pgfplotstableread{
1 0
5 969
10 806
20 719
40 656
80 622
}\datatableGammaMonoIter

\pgfplotstableread{
1 2893
5 2027
10 1483
20 1029
40 880
80 787
}\datatableGammaSplitIter

\pgfplotstableread{
1 0
5 1342.9570915699005
10 1554.938975572586
20 1303.4342246055603
40 1192.6504316329956
80 1157.1092596054077
}\datatableGammaMonoTime

\pgfplotstableread{
1 1129.1626975536346
5 731.7424337863922
10 713.6887679100037
20 473.86979150772095
40 387.60899901390076
80 364.61568999290466
}\datatableGammaSplitTime

\subsection{Dependence on swelling parameter}\label{sec:swell}
In addition, both the monolithic and splitting methods were evaluated using various swelling parameter values $\xi$. The findings mirror those observed with the surface tension parameter: an increase in the swelling parameter increases the coupling strength of the problem, subsequently elevating the computational demands for the linearization methods. Consistently, the splitting method demonstrated greater robustness in comparison to the monolithic method, successfully converging at higher coupling strengths. However, it is noteworthy that at elevated swelling parameter values ($\xi\geq 1.5$), instances were identified where the splitting method failed to converge as well. An unexplored potential approach, not covered in this paper, is the adoption of alternative time-discretization methods, as suggested in \cite{Storvik2023}.

\pgfplotstableread{
0 0
1 0
80 0
}\datatableEmpty

\pgfplotstableread{
0.01 608
0.1 623
0.25 969
0.5 0
0.75 0
1.0 0
10 0
}\datatableXiMonoIter

\pgfplotstableread{
0.01 603
0.1 732
0.25 2027
0.5 3088
0.75 4092
1.0 0
10 0
}\datatableXiSplitIter

\pgfplotstableread{
0.01 460.1577422618866
0.1 786.6613414287567
0.25 1240.788411617279
0.5 0
0.75 0
1.0 0
10 0
}\datatableXiMonoTime

\pgfplotstableread{
0.01 183.1707727909088
0.1 247.7530107498169
0.25 694.9803464412689
0.5 1115.081544160843
0.75 1480.4271116256714
1.0 0
10 0
}\datatableXiSplitTime

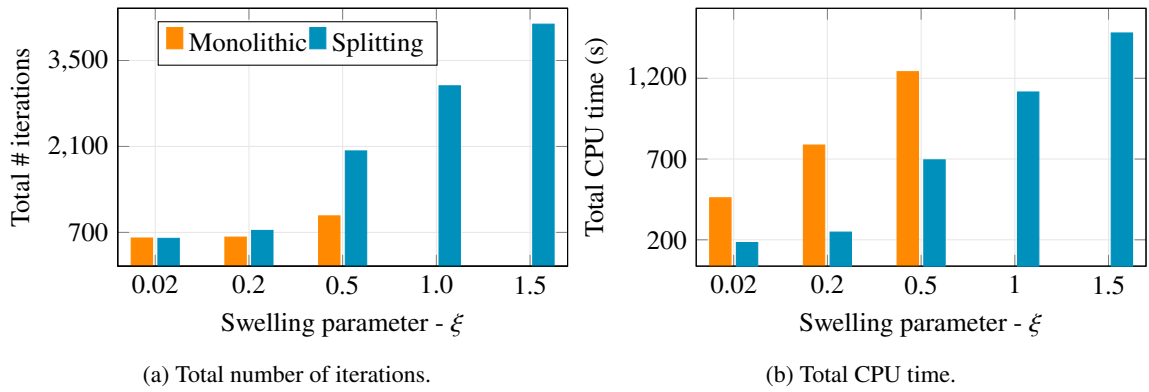
\begin{figure}
\begin{subfigure}{0.5\textwidth}
  \begin{tikzpicture}
   \pgfplotsset{ybar stacked, ymin=500.000, ymax=4000.000, width=\textwidth, enlargelimits=0.1,
   grid=both, grid style={line width=.2pt, draw=gray!20}, 
   xmin=0, xmax = 4,
 }
 
 \begin{axis}
 [
 ylabel={Total \# iterations}, 
 xticklabels={0.02, 0.2, 0.5, 1.0, 1.5, 2, 20},
 xlabel = {Swelling parameter - $\xi$},
 ytick={700, 2100, 3500}, 
 height=5cm,
 xtick=data,
 bar width=8pt,
 bar shift=-5pt, 
 height=5cm,
 scaled y ticks = false
 ]  
 
 \pgfplotsinvokeforeach {1}{
     \addplot[color=neworange1, fill=neworange1] table [x expr=\coordindex, y index=#1] {\datatableXiMonoIter};}
     
\end{axis} 

 \begin{axis}
 [
 bar width=8pt,
 bar shift=5pt, 
 xtick=\empty,
 ytick=\empty, 
 height=5cm
 ]  
 
 \pgfplotsinvokeforeach {1}{
     \addplot[color=newblue1, fill=newblue1] table [x expr=\coordindex, y index=#1] {\datatableXiSplitIter};}
     
\end{axis} 

 \begin{axis}
 [
 bar width=0pt,
 bar shift=10pt, 
 legend columns=2,
 legend style={at={(0.4,0.95)}, anchor=north}, 
 xtick=\empty,
 legend cell align=left,
 legend style={
     /tikz/column 2/.style={
         column sep=-1.3pt,
         row sep=-20pt 
     },
     /tikz/column 6/.style={ 
         column sep=-1.3pt,
         row sep=-20pt 
     }},
 ytick=\empty,
 height=5cm
 ]  
 \pgfplotsinvokeforeach {1}{
     \addplot[color=neworange1, fill=neworange1] table [x expr=\coordindex, y index=#1] {\datatableEmpty};}
 \pgfplotsinvokeforeach {1}{
     \addplot[color=newblue1, fill=newblue1] table [x expr=\coordindex, y index=#1] {\datatableEmpty};}

 \legend{
 Monolithic, 
 Splitting, 
 }
 \end{axis}

 \end{tikzpicture}
 
 \caption{Total number of iterations.}
 \label{fig:xiiterations}
 \end{subfigure}
 \begin{subfigure}{0.5\textwidth}
\begin{tikzpicture}
   \pgfplotsset{ybar stacked, ymin=170.000, ymax=1500.000, width=\textwidth, enlargelimits=0.1,
   grid=both, grid style={line width=.2pt, draw=gray!20}, 
   xmin=0, xmax = 4,
 }
 
 \begin{axis}
 [
 ylabel={Total CPU time (s)}, 
 xticklabels={0.02, 0.2, 0.5, 1, 1.5, 2},
 xlabel = {Swelling parameter - $\xi$},
 ytick={200, 700, 1200, 1700}, 
 height=5cm,
 xtick=data,
 bar width=8pt,
 bar shift=-5pt, 
 height=5cm,
 scaled y ticks = false
 ]  
 
 \pgfplotsinvokeforeach {1}{
     \addplot[color=neworange1, fill=neworange1] table [x expr=\coordindex, y index=#1] {\datatableXiMonoTime};}
     
\end{axis} 

 \begin{axis}
 [
 bar width=8pt,
 bar shift=5pt, 
 xtick=\empty,
 ytick=\empty, 
 height=5cm
 ]  
 
 \pgfplotsinvokeforeach {1}{
     \addplot[color=newblue1, fill=newblue1] table [x expr=\coordindex, y index=#1] {\datatableXiSplitTime};}
     
\end{axis} 



 \end{tikzpicture}
 \caption{Total CPU time.}
 \label{fig:xitime}
 \end{subfigure}
 \caption{The total number of iterations and CPU time to complete the simulation for different values of the swelling parameter $\xi$. The material parameters from Table~\ref{tab:1} are applied with surface tension parameter $\gamma = 5$.}
 \label{fig:xi}
\end{figure}

\section{Conclusion}
This study presents an investigation of both monolithic and splitting methods across varying conditions, specifically focusing on their performance under different surface tension and swelling parameter values. The results demonstrate a clear trend: as the coupling strength increases, the computational complexity for the linearization methods also rises. However, the splitting method consistently outperforms the monolithic approach in terms of robustness and efficiency. It converges more reliably across a broader range of parameter values, especially in scenarios with higher coupling strengths, although limitations are observed here as well. 

We note that this is a preliminary study, and there are several potential areas for further research. These include the exploration of alternative time-discretization methods for improved conditioning of the nonlinear systems, and potentially alternative ways to decouple the system of equations. In particular, different ways of stabilizing the splitting scheme with a hope of achieving even better robustness will be valuable. Additionally, a further mathematical study of the solution strategies will be of interest.

\subsection*{Acknowledgements}
Funded in part by Deutsche Forschungsgemeinschaft (DFG, German Research Foundation) under Germany's Excellence Strategy - EXC 2075 – 390740016. CR acknowledges the support by the Stuttgart Center for Simulation Science (SC SimTech). JWB acknowledges support from the UoB Akademia-project FracFlow. 
FAR acknowledges the support of the VISTA program, The Norwegian Academy of Science and Letters and Equinor.

\bibliographystyle{unsrt}
\bibliography{bibliography}

\end{document}